\documentclass[12pt, oneside]{article}
\usepackage{amsmath,amsfonts,amssymb,amscd,theorem}
\date{}
\newtheorem{thm}{Theorem}
\newtheorem{lem}{Lemma}
\newtheorem{defi}{Definition}

\newtheorem{prop}{Proposition}

\begin{document}

\title{A new proof of a theorem of Petersen}

\author{Yi-Hu Yang (SJTU)\footnote{Supported partially by NSF of China (No.11171253)}
~and Yi Zhang (Tongji)}

\maketitle

\begin{abstract}
Let $M$ be an $n$-dimensional complete Riemannian manifold with Ricci curvature $\ge n-1$. In \cite{colding1, colding2},
Tobias Colding, by developing some new techniques, proved that the following three condtions: 1) $d_{GH}(M, S^n)\to 0$;
2) the volume of $M$ ${\text{Vol}}(M)\to{\text{Vol}}(S^n)$; 3) the radius of $M$ ${\text{rad}}(M)\to\pi$ are equivalent.
In \cite{peter}, Peter Petersen, by developing a different technique, gave the 4-th equivalent condition,
namely he proved that the $n+1$-th eigenvalue of $M$ $\lambda_{n+1}(M)\to n$ is also equivalent to
the radius of $M$ ${\text{rad}}(M)\to\pi$, and hence the other two.
In this note, we give a new proof of Petersen's theorem by utilizing Colding's techniques.
\end{abstract}



\section{Introduction}
This note is a by-product of understanding \cite{colding1, colding2, peter}.
Let $M$ be an $n$-dimensional complete Riemannian manifold with Ricci curvature $\ge n-1$.
In \cite{colding1, colding2}, T. Colding proved that the following three conditions are equivalent:
1) $d_{GH}(M, S^n)\to 0$; 2) the volume of $M$ ${\text{Vol}}(M)\to{\text{Vol}}(S^n) $;
3) the radius of $M$ ${\text{rad}}(M)\to\pi$.
To this end, he developed some new techniques and got some local $L^2$-estimates of distances
and angles (for details, see \S 2). On the other hand, by developing a completely different technique,
Peter Petersen later got the 4-th equivalent condition \cite{peter}, i.e. proved the following theorem.

\vskip.5cm
\noindent
{\bf Petersen's Theorem:} {\it Let $M$ be an $n$ dimensional complete Riemannian manifold
with Ricci curvature $Ric_M\ge n-1$. Then the radius $rad(M)$ of $M$
is close to $\pi$ if and only if the $n+1$-th eigenvalue $\lambda_{n+1}(M)$ of $M$ is close to $n$.}

\vskip.5cm
The aim of this note is to give a new proof of Petersen's theorem by utilizing Colding's techniques.
Here, we also want to thank Peter Petersen for his kind comments.

\section{Colding's local $L^2$-estimates of distances and angles}
In this section, we recall Colding's integral estimates of distances and angles. We first fix some notation.
Let $N$ be a closed Riemannian manifold. For $p\in [1, \infty)$ and $f\in L^p(N)$, set
\[
\Vert f\Vert_p=\big({\frac{1}{{\text{Vol}}(N)}}\int_N|f|^p{\text d}vol\big)^{\frac 1 p}.
\]
In the sequel, we will always omit "${\text d}vol$" in the integral. Let $f$ be a Lipschitz function and $\nabla f$
the gradient of $f$, then by $\Vert f\Vert_{2,1}$ denote the $(2,1)$-Sobolev norm of $f$, i.e.
\[
\Vert f\Vert_{(2,1)}=\big({\frac{1}{{\text{Vol}}(N)}}\int_N|f|^2+{\frac{1}{{\text{Vol}}(N)}}\int_N|\nabla f|^2\big)^{\frac 1 2}.
\]
For $f\in C^2(N)$, define the Hessian of $f$ as follows
\[
{\text{Hess}}(f)(u,v)=<\nabla_u\nabla f, v>,
\]
and the laplacian as follows
\[
\Delta f={\text{Tr}}{\text{Hess}}(f).
\]
So, $\Delta$ is negative semi-definite.

\vskip.5cm
Recall also that $f$ is an eigenfunction with eigenvalue $\lambda$
if $\Delta f+\lambda f=0$. Let $SN$ be the unit tangent bundle of $N$, $g^t$ is the geodesic flow and
$\pi:SN\to N$ the corresponding projection. For $A\subset N$, set $SA=\pi^{-1}(N)$. For $SN$, one has
an obvious Riemannian measure which is induced from that of $N$. So, we will integrate on various
sets, e.g. $N, SN, SA\times [0, l]$, etc.

\vskip.5cm
If $V$ is a finite dimensional vector space with an inner product $<\cdot, \cdot>$ and $\omega$ is a
bilinear from on $V$, then we set
\[
|\omega |^2=\Sigma_{i,j}(\omega(e_i, e_j))^2,
\]
where $\{e_i\}$ is an orthonormal basis of $V$.

\vskip.5cm
In this note, unless stated otherwise, we always assume that $M$ is an $n$-dimensional complete Riemannian
manifold with the Ricci curvature satisfying $Ric_M\ge n-1, n\ge 2$.

\vskip.5cm
For $f\in C^{\infty}(M)$ and $l>0$, define a function $h_f\in C^{\infty}(SM\times [0, l])$ as follows
\begin{eqnarray}
h_f(v, t)=f(\gamma_v(0))\cos t+{\frac{f(\gamma_v(l))-f(\gamma_v(0))\cos l}{\sin l}}\sin t,
\end{eqnarray}
where $\gamma_v(t)$ is the geodesic at $\pi(v)$ along the direction $v$. The function obviously satisfies
\begin{eqnarray}
{\frac{\partial^2h_f}{\partial t^2}}=-h_f
\end{eqnarray}
and is also determined uniquely under the conditions: $h_f(v, 0)=f(\gamma_v(0))$, $h_f(v, l)=f(\gamma_v(l))$.
We remark that if ${\text{Hess}}(f)=-fg$ (here $g$ is the Riemannian metric on $M$), then
$h_f(v, t)=f(\gamma_v(t))$.
\begin{prop} (\cite{colding1}, Lemma 1.4)
Given $\epsilon, k$, and $l\in [{\frac{\pi}{2}}, \pi]$, there exists a positive number
$\delta=\delta(\epsilon, k, l, n)$ such that if $f\in C^{\infty}(M)$ with $\Vert f\Vert_2\le k$
and $\Vert \Delta f+nf\Vert_2\le\delta$, then
\begin{eqnarray}
{\frac{1}{l{\text{Vol}}(SM)}}\int_{SM}\int_0^l|f(\gamma_v(t))-h_f(v, t)|^2dt<\epsilon,
\end{eqnarray}
and
\begin{eqnarray}
{\frac{1}{l{\text{Vol}}(SM)}}\int_{SM}\int_0^l\big|{\frac{\partial f(\gamma_v(t))}{\partial t}}
-{\frac{\partial h_f(v, t)}{\partial t}}\big|^2dt<\epsilon.
\end{eqnarray}
\end{prop}

\vskip.5cm
In the sequel of this note, we will use $\psi(\delta |\cdot,\cdots,\cdot)$ to denote
a certain positive function depending on $\delta$ and some additional parameters
such that when these parameters are fixed, $\lim_{\delta\to 0}\psi(\delta |\cdot,\cdots,\cdot)=0$.

\vskip.5cm
For any $f\in C^{\infty}(M)$, we define another function $g_f$ on $SM\times [0, l]$ as follows
\begin{eqnarray}
g_f(v, t)=<\nabla f, v>\sin t+f(\pi(v))\cos t.
\end{eqnarray}
Similar to the $h_f$ before, it is easily to see that $g_f(v, t)$ is uniquely determined by
\begin{eqnarray}
&&{\frac{\partial^2g_f}{\partial t^2}}=-g_f,\\
&&g_f(v, 0)=f(\pi(v)),\\
&&{\frac{\partial g_f}{\partial t}}(v, 0)=<\nabla f, v>.
\end{eqnarray}

The following proposition also is essentially due to Colding. We only make a little bit modification so that
we can apply it conveniently in the proof of the next section.
\begin{prop} (\cite{colding2}, Proposition 4.5)
Given $\epsilon, k, \bar{A}, C>0$ and $l\in [{\frac \pi 2}, \pi)$, there exists a positive
number $\delta=\delta(n, k, \epsilon, l, {\bar A}, C)$ such that
if $A\subset M$ with ${\frac{Vol(A)}{Vol(M)}}\ge{\bar A}$ and $f\in C^{\infty}(M)$ with $\Vert f\Vert_2\le k$,
$\Vert \Delta f+nf\Vert_2\le\delta$, $\max|f|$, and $\max|\nabla f|\le C$, then we have
\begin{eqnarray}
{\frac{1}{l^3{\text{Vol}}(SA)}}\int_{SA}\int_0^l|g_f(v,t)-f(\gamma_v(t))|^2dt<\epsilon,
\end{eqnarray}
and
\begin{eqnarray}
{\frac{1}{l{\text{Vol}}(SA)}}\int_{SA}\int_0^l\big|{{\frac{\partial g_f(v, t)}{\partial t}}
-\frac{\partial f(\gamma_v(t))}{\partial t}}\big|^2dt<\epsilon.
\end{eqnarray}
\end{prop}

\vskip.5cm
\noindent
{\it Proof.} The first estimate easily follows from the second one by integration.
By the definition of $h_f$ and the boundedness of $f$ and $|\nabla f|$,
for $\forall s\in [0, \epsilon l]$ we have
\begin{eqnarray*}
|h_f(v, s)-f(\gamma_v(s))|^2&\le&|h_f(v, s)-h_f(v, 0)|^2+|f(\gamma_v(s))-f(\gamma_v(0))|^2\\
&<&\psi(\epsilon|n, C).
\end{eqnarray*}
By Proposition 1, we also have, for some $\delta=\delta(\epsilon^2, k, l, n)$,
\[
{\frac{1}{l{\text{Vol}}(SM)}}\int_{\pi^{-1}T_{\epsilon l}(A)}
\int_0^{\epsilon l}\big|{{\frac{\partial h_f(v, t)}{\partial t}}
-\frac{\partial f(\gamma_v(t))}{\partial t}}\big|^2dt<\epsilon^2,
\]
where $\pi^{-1}T_{\epsilon l}(A)$ is the $\epsilon$-neighborhood of $A$. In particular,
$\exists s'\in[0, \epsilon l]$ such that
\[
{\frac{1}{{\text{Vol}}(SM)}}\int_{\pi^{-1}T_{\epsilon l}(A)}
\big|{{\frac{\partial h_f(v, s')}{\partial t}}
-\frac{\partial f(\gamma_v(s'))}{\partial t}}\big|^2<\epsilon.
\]
On the other hand, by the definition of $h_f$ and $g_f$, we have
\[
{\frac{\partial f(\gamma_v)}{\partial t}}(s')={\frac{\partial g_f(g^{s'}v, 0)}{\partial t}}.
\]
Also, for the boundedness of $f$ and $|\nabla f|$, we have
\[
\big|{\frac{\partial h_f(g^{s'}v, 0)}{\partial t}}-{\frac{\partial h_f(v, s')}{\partial t}}\big|<
\psi(\epsilon|n, \cdots).
\]
Thus, we have
\[
\big|{\frac{\partial h_f(g^{s'}v, 0)}{\partial t}}-{\frac{\partial g_f(g^{s'}v, 0)}{\partial t}}\big|<
\psi(\epsilon|n, \cdots).
\]
Then, the fact that $h_f(g^{s'}v£¬t)$ and $g_f(g^{s'}v, t)$ satisfy the same equation $f''=-f$,
together with $h_f(g^{s'}v, 0)=g_f(g^{s'}v, 0)$ implies
\[
{\frac{1}{{l\text{Vol}}(SM)}}\int_{\pi^{-1}T_{s'}(A)}\int_0^l
\big|{{\frac{\partial g_f(g^{s'}v, t)}{\partial t}}
-\frac{\partial h_f(g^{s'}v, t)}{\partial t}}\big|^2dt<\psi(\epsilon|n, C...).
\]
Using the fact that the geodesic flow is volume-preserving, we get
\[
{\frac{1}{{l\text{Vol}}(SM)}}\int_{SA}\int_0^l
\big|{{\frac{\partial g_f(v, t)}{\partial t}}
-\frac{\partial h_f(v, t)}{\partial t}}\big|^2dt<\psi(\epsilon|n, C...).
\]
Then, we have
\begin{eqnarray*}
{\frac{1}{{l\text{Vol}}(SA)}}\int_{SA}\int_0^l
\big|{{\frac{\partial g_f(v, t)}{\partial t}}
-\frac{\partial h_f(v, t)}{\partial t}}\big|^2dt&<&{\frac{{\text{Vol}}(SM)}{{\text{Vol}}(SA)}}
\psi(\epsilon|n, C...)\\
&=&\psi(\epsilon|n, C, \bar A...).
\end{eqnarray*}
Combining this with Proposition 1, we get the required inequality with $\epsilon$
replaced by $\psi(\epsilon|k, l, n, C, \bar A)$.

\vskip.5cm
\noindent
{\bf Remark}. When the set $A$ in the proposition is sufficiently small and $l$ sufficiently close to $\pi$,
$\gamma_v(t)$ can run over $M$, but $g_f(v, t)$ depends only on $f_{|A}$. So, such an $f$ has
much more geometric information of $M$. This is the key point of the proposition

\section{The proof of Petersen's theorem}
\subsection{$\lambda_{n+1}\to n$ implies ${\text{rad}}(M)\to\pi$}
We first give some preliminaries.
Let $p\in M$. Set $r(p)=\min\{r | B_r(p)\supset M\}$, the radius at $p$ of $M$.
Then the radius of $M$ is defined as
$$
{\text{rad}(M)}=\max_{p\in M}r(p).
$$
We also give the following
\begin{defi}
For $p, q\in M$ and $\delta, s>0$, set
\[
C_s=\{v\in SM | \pi(v)\in B_{\delta}(p), \exp (sv)\in B_{\delta}(q), d(\exp(sv), \pi(v))=s\}.
\]
\end{defi}
The following lemma is also due to Colding (\cite{colding1}, Lemma 2.3);
for convenience in the following proofs, we here give it a detailed proof.
\begin{lem}
$\forall\delta>0$ and $p, q\in M$, then there exists an $s'>0$ with
\[
{\text{Vol}}(C_{s'})\ge {\frac{n}{\pi^n}}\big(\frac{V_n(\delta){\text{Vol}}(M)}{\omega^n}\big)^2,
\]
where $V_n(\delta)$ denotes the volume of the ball with radius $\delta$ in the standard $n$-sphere $S^n$,
$\omega_n$ is the volume of $S^n$.
\end{lem}
{\it Proof}. Set
\[
C=\{v\in TM | \pi(v)\in B_{\delta}(p), \exp (v)\in B_{\delta}(q), d(\exp(v), \pi(v))=|v|\}.
\]
For $x\in M$, by Bishop's volume comparison, one has that the exponential map $\exp_x$ at $x$,
when restricted to the set
\[
C_x=\{v\in T_xM | d(\exp_x(v), x)=|v|\},
\]
is volume non-increasing. In particular, for all $x\in B_{\delta}(p)$, the Bishop-Gromov comparison theorem implies
\[
{\text{Vol}}_{\mathbb{R}^n}(T_xM\cap C)\ge {\text{Vol}}(B_\delta(q))\ge {\frac{V_n(\delta){\text{Vol}}(M)}{\omega_n}}.
\]
Thus, one has
\[
{\text{Vol}}(C)\ge \big({\frac{V_n(\delta){\text{Vol}}(M)}{\omega_n}}\big)^2.
\]
On the other hand, one has
\[
{\text{Vol}}(C)=\int_0^\pi{\text{Vol}}(C_s)s^{n-1}ds\le {\frac{\pi^n}{n}}\max_{s>0}{\text{Vol}}_SM(C_s).
\]
Therefore, there exists $s>0$ with
\[
{\text{Vol}}(C_s)\ge {\frac{n}{\pi^n}}\big(\frac{V_n(\delta){\text{Vol}}(M}{\omega^n}\big)^2.
\]
The lemma is obtained.

\vskip.5cm
Let $n\le\lambda_1\le\lambda_2\le\cdots\le\lambda_{n+1}$ be the first $n+1$ nonzero eigenvalues of
$M$ and $f_1, f_2, \cdots, f_{n+1}$ be the corresponding eigenfunctions respectively, i.e.
\[
\Delta f_i+\lambda_if_i=0, i=1, 2, \cdots, n+1;
\]
furthermore, we can assume that for $i\neq j$, $\int_Mf_if_j=0$.
We can normalize each $f_i$ to make it satisfy
\[
{\frac{1}{{\text{Vol}}(M)}}\int_Mf_i^2=1.
\]
Let $a_1, a_2, \cdots, a_{n+1}$ be $n+1$ real numbers satisfying $a_1^2+a_2^2+\cdots +a_{n+1}^2=1$.
Set $f=a_1f_1+a_2f_2+\cdots +a_{n+1}f_{n+1}$. Write $\lambda_{n+1}=n+\delta_1$. Since we only concern
the case of $\lambda_{n+1}$ near $n$, so WLOG, we can assume $\delta_1\le{\frac 1{n+1}}$.
From the definition of $f$, we have
\[
{\frac{1}{{\text{Vol}}(M)}}\int_Mf^2=1,
\]
and
\[
\Delta f+nf=\Sigma_{i=1}^{n+1}a_i(\Delta f_i+nf_i)=\Sigma_{i=1}^{n+1}a_i(n-\lambda_i)f_i.
\]
Furthermore, we have
\[
\Vert\Delta f+nf\Vert_2\le\delta_1.
\]
So, for all $a_1, a_2, \cdots, a_{n+1}$ with $a_1^2+\cdots+a_{n+1}^2=1$, we have
\[
\Vert\Delta f+nf\Vert_2\to 0 ~{\text{uniformly}},~ {\text{as}} ~\lambda_{n+1}\to n.
\]
Moreover, by the standard estimate of PDE of elliptic type, there exists a positive
constant $C$ (independent of $a_1, \cdots, a_{n+1}$) such that for the above $f$ one has
\[
\max_M|f|\le C,
\]
\[
\max_M|\nabla f|\le C,
\]
and
\[
\max_M|{\text{Hess}}(f)|\le C.
\]
In the following, when we mention a smooth $f$ on $M$,
we always mean such a function unless stated otherwise.

\vskip.5cm
Using these preliminaries and Colding's local integral estimates of distances and angels (Proposition 2),
we can now prove the following two lemmas.

\begin{lem}
For arbitrarily given $\epsilon$, there exists
$\delta'_1=\delta'_1(\epsilon)$, $\delta'_2=\delta'(\epsilon)$ and $\delta_3=\delta_3(\epsilon)$
such that if $\lambda_{n+1}\le n+\delta'_1$, and for some $p\in M$ and
some $\alpha\le{\frac{\epsilon}{8}}$, some smooth function $f$
on $M$ as mentioned before satisfies
\[
{\frac{1}{{\text{Vol}}(B_\alpha(p))}}\int_{B_\alpha(p)}|\nabla f|^2<\delta'_2
\]
and
\[
{\frac{1}{{\text{Vol}}(B_\alpha(p))}}\int_{B_\alpha(p)}|f|^2<\delta_3,
\]
then $r(p)>\pi-\epsilon$.
\end{lem}
{\it Proof}: Assume $r(p)\le\pi-\epsilon$. Then $B_{\pi-\epsilon}(p)\supset M$.
Applying Proposition 2 with $A=B_{\alpha}(p)$ and $l=\pi-{}\frac{\epsilon}{2}$,
we can choose a $\delta'_1=\delta_1'(\epsilon)=\delta(n, {\frac{\epsilon}{2(\pi-{\frac{\epsilon}{2}})}},
\pi-{\frac{\epsilon}{2}}, {\frac{{\text{Vol}}(S^n(\alpha))}{\omega_n}}, C)$ satisfying
\begin{eqnarray*}
{\frac{1}{(\pi-{\frac{\epsilon}{2}}){\text{Vol}}(SB_\alpha(p))}}
\int_{SB_\alpha(p)}\int_0^{\pi-{\frac{\epsilon}{2}}}|g_f(v, t)-f(\gamma_v(t))|^2\le{\frac{\epsilon}{2}}.
\end{eqnarray*}
On the other hand, we have
\begin{eqnarray*}
1&=&{\frac{1}{{\text{Vol}}(M)}}\int_Mf^2\\
&\le&{\frac{(\pi-{\frac{\epsilon}{2}})^{n-1}}{{\text{Vol}}(M){\text{Vol}}(B_{\alpha}(p))}}
\int_{SB_\alpha(p)}\int_0^{\pi-{\frac{\epsilon}{2}}}|f(\gamma_v(t))|^2\\
&\le&{\frac{C'}{(\pi-{\frac{\epsilon}{2}}){\text{Vol}}(SB_\alpha(p))}}
\int_{SB_\alpha(p)}\int_0^{\pi-{\frac{\epsilon}{2}}}\big(|f(\gamma_v(t))-g_f(v, t)|^2+|g_f(v, t)|^2\big)\\
&\le&C''({\frac \epsilon 2}+\delta_2+\delta_3).
\end{eqnarray*}
So, when $\epsilon$, $\delta_2$, and $\delta_3$ are sufficiently small, this derives a contradiction. The lemma is obtained.

\begin{lem}
For arbitrarily given $\epsilon$ and $\delta_3$, set $\alpha_0=\min\{{\frac{C'\delta_3\epsilon^2}{80C}}, {\frac{\epsilon}{8}}\}$
(we take $C'\le{\frac{1}{10}}$ later). Then there exist $\delta''_1=\delta''_1(\epsilon, \alpha_0, \delta_3)$ and
$\delta''_2=\delta''_2(\epsilon, \delta_3)$ satisfying that if $\lambda_{n+1}\le n+\delta''_1$ and, for some smooth
function $f$ as mentioned before, $p\in M$ and some $\alpha\le\alpha_0$,
\[
{\frac{1}{{\text{Vol}}(B_\alpha(p))}}\int_{B_\alpha(p)}|\nabla f|^2<\delta''_2
\]
and
\[
{\frac{1}{{\text{Vol}}(B_\alpha(p))}}\int_{B_\alpha(p)}|f|^2\ge\delta_3,
\]
then $r(p)>\pi-\epsilon$.
\end{lem}
{\it Proof.} Assume $r(p)\le\pi-\epsilon$. Let $f$ be a smooth function mentioned before, $p\in M$ and
$\alpha\le\alpha_0$. For any $v, v'\in SM$, $t,t'\in[0, \pi-{\frac{\epsilon}{4}}]$,
we have
\begin{eqnarray*}
&&|f(\gamma_v(t))-g_f(v, t)|\\
&=&|f(\gamma_{v'}(t'))-g_f(v', t')-g_f(v, t)+g_f(v', t')-f(\gamma_{v'}(t'))+f(\gamma_v(t))|\\
&\ge&|g_f(v, t)-g_f(v', t')|-|f(\gamma_{v'}(t'))-g_f(v', t')|-|f(\gamma_{v'}(t'))-f(\gamma_v(t))|,
\end{eqnarray*}
so, we have
\begin{eqnarray}
&&|f(\gamma_v(t))-g_f(v, t)|+|f(\gamma_{v'}(t'))-g_f(v', t')|  \nonumber\\
&\ge&|g_f(v, t)-g_f(v', t')|-|f(\gamma_{v'}(t'))-f(\gamma_v(t))|.
\end{eqnarray}

Since ${\frac{1}{{\text{Vol}}(B_\alpha(p))}}\int_{B_\alpha(p)}|f|^2\ge\delta_3$, there exists
a $q_1\in B_\alpha(p)$ such that $|f(q_1)|\ge\delta_3$. From this and the fact that $|\nabla f|\le C$,
$\exists r_0=r_0(\epsilon, \delta_3, C)\le{\frac{\epsilon}{32}}$ such that for any $r\le r_0$ and $q_2\in B_{r}(q_1)$,
$|f(q_2)|\ge {\frac{\delta_3}{2}}$.

\vskip.5cm
We first derive a lower bound of the term $|g_f(v, t)-g_f(v', t')|$
for some $v, v', t, t'$.
It is clear that, for any $v\in SB_r(q_1)$, $|f(\pi(v'))|\ge{\frac{\delta_3}{4}}$
for any $v'\in SB_\beta(\pi(v))$ as $\beta$ is sufficiently small ($\beta$ will be fixed in the following).
Combining these with $|\nabla f|\le C$ and $|{\text{Hess}}(f)|\le C$, we have
\begin{eqnarray*}
&&|g_f(v, t)-g_f(v', t')|\nonumber\\
&\ge&|f(\pi(v))||\cos t-\cos t'|-|f(\pi(v))-f(\pi(v'))|-|\nabla f(\pi(v))|-|\nabla f(\pi(v'))|\nonumber\\
&\ge&|f(\pi(v))||\cos t-\cos t'|-2C\beta-2|\nabla f(\pi(v))|.
\end{eqnarray*}
Taking $t\in [\pi-{\frac{\epsilon}{2}}, \pi-{\frac{\epsilon}{4}}]$ and $t'\in [0, \pi-{\frac{3\epsilon}{4}}]$,
we then have
\begin{eqnarray*}
|g_f(v, t)-g_f(v', t')|\ge C'\delta_3\epsilon^2-2C\beta-2|\nabla f(\pi(v))|,
\end{eqnarray*}
here $C'\le {\frac{1}{10}}$.
Take $\beta\le\min\{{\frac{C'\delta_3\epsilon^2}{20C}}, {\frac{\epsilon}{32}}\}$. We then have
\begin{eqnarray*}
|g_f(v, t)-g_f(v', t')|\ge {\frac{9}{10}}C'\delta_3\epsilon^2-2|\nabla f(\pi(v))|.
\end{eqnarray*}
On the other hand, since ${\frac{1}{{\text{Vol}}(B_\alpha(p))}}\int_{B_\alpha(p)}|\nabla f|^2\le\delta''_2$
($\delta''_2$ will be fixed in the following), $|{\text{Hess}}(f)|\le C$ and $v\in SB_{r}(q_1)$, we have
\[
|\nabla f(\pi(v))|\le Cr+2C\alpha+\delta''_2.
\]
Taking $\delta''_2\le{\frac{C'\delta_3\epsilon^2}{80}}$
and $r_0\le\min\{{\frac{C'\delta_3\epsilon^2}{80C}}, {\frac{\epsilon}{32}}\}$,
we then have
\[
|\nabla f(\pi(v))|\le{\frac{C'\delta_3\epsilon^2}{20}}.
\]
Consequently, we have
\begin{eqnarray}
|g_f(v, t)-g_f(v', t')|\ge {\frac{8}{10}}C'\delta_3\epsilon^2,
\end{eqnarray}
for any $(v, t)\in SB_r(q_1)\times[\pi-{\frac{\epsilon}{2}}, \pi-{\frac{\epsilon}{4}}]$
and any $(v', t')\in SB_\beta(\pi(v))\times[0, \pi-{\frac{3\epsilon}{4}}]$.

\vskip.5cm
Next, we want to derive an upper bound of the term $|f(\gamma_v(t))-g_f(v, t)|$
for some suitable $(v, t)\in SB_r(q_1)\times[\pi-{\frac{\epsilon}{2}}, \pi-{\frac{\epsilon}{4}}]$.
Setting $l=\pi-{\frac{\epsilon}{4}}$ and applying Proposition 2 on $A=B_{\frac{\epsilon}{4}}(p)$, we then have
\begin{eqnarray*}
{\frac{1}{(\pi-{\frac{\epsilon}{4}})^3{\text{Vol}}(SB_{\frac{\epsilon}{4}}(p))}}
\int_{SB_{\frac{\epsilon}{4}}(p)}\int_0^{\pi-{\frac{\epsilon}{4}}}|g_f(v,t)-f(\gamma_v(t))|^2dt<\psi_0(\delta''_1 | \epsilon),
\end{eqnarray*}
here $\psi_0(\delta''_1|\epsilon)\to 0$ as $\delta''_1\to 0$ ($\delta''_1$ will be fixed in the following).

We also remark that $\alpha\le\alpha_0\le{\frac{\epsilon}{8}}$ and $r\le r_0\le{\frac{\epsilon}{32}}$,
and hence $SB_r(q_1)\subset SB_{\frac{\epsilon}{4}}(p)$.
So, by the volume comparison theorem, we have
\begin{eqnarray}
{\frac{1}{\epsilon{\text{Vol}}(SB_{r}(q_1))}}
\int_{SB_{r}(q_1)}\int_{\pi-{\frac{\epsilon}{2}}}^{\pi-{\frac{\epsilon}{4}}}
|g_f(v,t)-f(\gamma_v(t))|^2dt<\psi_1(\delta''_1|\epsilon).
\end{eqnarray}
Thus, for sufficiently small $\delta''_1$ and some
$(v, t)\in SB_r(q_1)\times [\pi-{\frac{\epsilon}{2}}, \pi-{\frac{\epsilon}{4}}]$, we have
\begin{eqnarray}
|f(\gamma_v(t))-g_f(v,t)|<\sqrt{\psi_1(\delta''_1|\epsilon)}\le {\frac{1}{10}}C'\delta_3\epsilon^2.
\end{eqnarray}
Combining this with (11) and (12), we have, for sufficiently small $\delta''_1$,
\begin{eqnarray}
|f(\gamma_{v'}(t'))-g_f(v', t')|+|f(\gamma_{v'}(t'))-f(\gamma_v(t))|\ge {\frac{7}{10}}C'\delta_3\epsilon^2,
\end{eqnarray}
for some $(v, t)\in SB_r(q_1)\times[\pi-{\frac{\epsilon}{2}}, \pi-{\frac{\epsilon}{4}}]$
and any $(v', t')\in SB_\beta(\pi(v))\times[0, \pi-{\frac{3\epsilon}{4}}]$.

\vskip.5cm
Now, we discuss the terms $|f(\gamma_{v'}(t'))-g_f(v', t')|$ and
$|f(\gamma_{v'}(t'))-f(\gamma_v(t))|$ in (15).
To this end, we need Lemma 1. For convenience, we first set
\begin{eqnarray*}
C(v, t)&=&\{(\bar v, s)\in SM\times[0, \pi]~|~\pi(\bar v)\in B_\beta(\pi(v)), \exp(s\bar v)\in B_\beta(\gamma_v(t)),\\
&&{\text{and}}~d(\exp(s\bar v), \pi(\bar v))=s\};\\
C'(v, t)&=&\{\bar v\in TM~|~\pi(\bar v)\in B_\beta(\pi(v)), \exp(\bar v)\in B_\beta(\gamma_v(t)),\\
&&{\text{and}}~d(\exp(\bar v), \pi(\bar v))=|\bar v|\};\\
C(s; v, t)&=&\{\bar v\in SM~|~\pi(\bar v)\in B_\beta(\pi(v)), \exp(s\bar v)\in B_\beta(\gamma_v(t)),\\
&&{\text{and}}~d(\exp(s\bar v), \pi(\bar v))=s\}.
\end{eqnarray*}

By the definition of $C(v, t)$ and the assumption $r(p)\le\pi-\epsilon$, we have, if $(\bar v, s)\in C(v, t)$,
$$
s\le\pi-\epsilon+\alpha+r+\beta+\beta\le\pi-{\frac{25\epsilon}{32}}<\pi-{\frac{3\epsilon}{4}}.
$$
So, $C(v, t)\subset SB_\beta(\pi(v))\times[0, \pi-{\frac{3\epsilon}{4}}]
\subset SB_{\frac{\epsilon}{4}}(p)\times[0, \pi-{\frac{3\epsilon}{4}}]$.
Similarly, for $\bar v\in C'(v, t)$, $|\bar v|\le\pi-{\frac{3\epsilon}{4}}$;
and for $s\ge\pi-{\frac{3\epsilon}{4}}$, $C(s; v, t)=\phi$. We also have
\[
{\text{Vol}}(C(v, t))=\int_0^\pi{\text{Vol}}(C(s; v, t))ds
=\int_0^{\pi-{\frac{3\epsilon}{4}}}{\text{Vol}}(C(s; v, t))ds
\]
and
\[
{\text{Vol}}(C'(v, t))=\int_0^\pi{\text{Vol}}(C(s; v, t))s^{n-1}ds
=\int_0^{\pi-{\frac{3\epsilon}{4}}}{\text{Vol}}(C(s; v, t))s^{n-1}ds.
\]

On the other hand, by the proof of Lemma 1, we have
\[
{\text{Vol}}(C'(v, t))\ge \big({\frac{V_n(\beta){\text{Vol}}(M)}{\omega_n}}\big)^2.
\]
So, we have
\begin{eqnarray*}
{\text{Vol}}(C(v, t))&=&\int_0^{\pi-{\frac{3\epsilon}{4}}}{\text{Vol}}(C(s; v, t))ds\\
&\ge&{\frac{1}{(\pi-{\frac{3\epsilon}{4}})^{n-1}}}\int_0^{\pi-{\frac{3\epsilon}{4}}}{\text{Vol}}(C(s; v, t))s^{n-1}ds\\
&=&{\frac{1}{(\pi-{\frac{3\epsilon}{4}})^{n-1}}}{\text{Vol}}(C'(v, t))\\
&\ge&{\frac{1}{(\pi-{\frac{3\epsilon}{4}})^{n-1}}}\big({\frac{V_n(\beta){\text{Vol}}(M)}{\omega_n}}\big)^2.
\end{eqnarray*}
Again applying Proposition 2 with $A=M$ and $l=\pi-{\frac{\epsilon}{4}}$, we have
\begin{eqnarray*}
{\frac{1}{{\text{Vol}}(SM)}}
\int_{SM}\int_{0}^{\pi-{\frac{\epsilon}{4}}}
|g_f(\bar v,s)-f(\gamma_{\bar v}(s))|^2ds<\psi_2(\delta''_1|\epsilon).
\end{eqnarray*}
Since $C(v, t)\subset SB_\beta(\pi(v))\times[0, {\pi-\frac{3\epsilon}{4}}]
\subset SM\times[0, {\pi-\frac{\epsilon}{4}}]$, we have
\begin{eqnarray*}
{\frac{1}{{\text{Vol}}(C(v, t))}}
\int_{0}^{\pi-{\frac{3\epsilon}{4}}}\big(\int_{C(s; v, t)}
|g_f(\bar v,s)-f(\gamma_{\bar v}(s))|^2\big)ds<\psi_3(\delta''_1|\epsilon).
\end{eqnarray*}
So, for sufficiently small $\delta''_1$ (which now can be fixed), there exists an $(v', t')\in C(v, t)
\subset SB_\beta(\pi(v))\times[0, {\pi-\frac{3\epsilon}{4}}]$ satisfying
\begin{eqnarray}
|g_f(v', t')-f(\gamma_{v'}(t'))|^2<\sqrt{\psi_3(\delta''_1|\epsilon)}\le {\frac{1}{10}}C'\delta_3\epsilon^2.
\end{eqnarray}
Since $(v', t')\in C(v, t)$, (by the definition of $C(v, t)$) we have $\gamma_{v'}(t')\in B_\beta(\gamma_v(t))$, so
\begin{eqnarray}
|f(\gamma_{v'}(t'))-f(\gamma_v(t))|\le C\beta\le{\frac{1}{10}}C'\delta_3\epsilon^2.
\end{eqnarray}
Combining (17) with (15) and (16), we derive a contradiction. The lemma is obtained.

\vskip.5cm
Clearly, Lemma 2 and Lemma 3 imply the following
\begin{thm}
For arbitrarily given $\epsilon>0$, there exist $\delta_1=\delta_1(\epsilon)>0$,
and $\delta_2=\delta_2(\epsilon)>0$
such that if $\lambda_{n+1}\le n+\delta_1$, and for some $p\in M$, some
smooth function $f$ mentioned before, and some $\alpha=\alpha(\epsilon)>0$
\[
{\frac{1}{{\text{Vol}}(B_\alpha(p))}}\int_{B_\alpha(p)}|\nabla f|^2<\delta_2,
\]
then $r(p)>\pi-\epsilon$.
\end{thm}
{\it Proof of 3.1}. Take $\epsilon>0$ (sufficiently small) and $p\in M$ arbitrarily. Consider the gradient vector
$\nabla f_1, \nabla f_2, \cdots, \nabla f_{n+1}$ at $p$ of the eigenfunctions of $M$.
There exist $n+1$ real numbers $a_1, a_2, \cdots, a_{n+1}$ with $\sum_{i=1}^{n+1}a_i^2=1$
satisfying $\sum_{i=1}^{n+1}a_i\nabla f_i=0$ at $p$. Set $f=\sum_{i=1}^{n+1}a_i f_i$. So, $(\nabla f)(p)=0$.
Thus, for the $\delta_2$ in Theorem 1, we have, for $\alpha\le{\frac{\sqrt\delta_2}{C}}$,
\[
|\nabla f|<\sqrt\delta_2, ~~{\text{on}}~B_\alpha(p).
\]
So, we have
\[
{\frac{1}{{\text{Vol}}(B_\alpha(p))}}\int_{B_\alpha(p)}|\nabla f|^2<\delta_2.
\]
We remark that the choose of $\delta_2$ and $\alpha$ are independent of $p$.
So, by Theorem 1, as $\lambda_{n+1}$ is sufficiently close to $n$, $r(p)>\pi-\epsilon$,
and hence ${\text{rad}}(M)>\pi-\epsilon$. The proof is finished.

\subsection{${\text{rad}}(M)\to\pi$ implies $\lambda_{n+1}\to n$}
Colding's theorem \cite{colding1} says that {${\text{rad}}(M)\to\pi$  is
equivalent to $d_{GH}(M, S^n)\to 0$. So, we just need to prove that $d_{GH}(M, S^n)\to 0$ implies $\lambda_{n+1}\to n$.
To do this, we need the following result of Colding (\cite{colding1}, Lemma 1.10).

\begin{lem}
$\forall$ $\epsilon>0$, $\exists$ $\delta=\delta(\epsilon, n)>0$, such that if there exist some $p, q\in M$
with $d(p, q)>\pi-\delta$, then there exists an $f\in C^\infty (M)$ with $\Vert f\Vert_2\le 1$,
$\Vert\Delta f+nf\Vert_2<\epsilon$, and $\Vert f-g\Vert_{2,1}<\epsilon$, here $g(x)=\cos d(p, x)$.
\end{lem}

\noindent
{\bf Remark.} Actually, for the $f$ in the above lemma, we can further assume $\int_Mf=0$.
In fact, for the above $f$, we have
\begin{eqnarray*}
\big|{\frac{1}{{\text{Vol}}(M)}}\int_Mf\big|&\le&{\frac 1 n}\big|{\frac{1}{{\text{Vol}}(M)}}\int_M(\triangle f+nf)\big|\\
&\le&{\frac{1}{n}}\Vert\Delta f+nf\Vert_2<{\frac \epsilon n}.
\end{eqnarray*}
Set $\bar f=f-{\frac{1}{{\text{Vol}}(M)}}\int_Mf$ with $\int_M\bar f=0$. Then $\Vert \bar f\Vert_2\le 1+\epsilon$,
$\Vert\Delta{\bar f}+n{\bar f}\Vert_2<2\epsilon$, and $\Vert \bar f-g\Vert_{2,1}<2\epsilon$.
So, we can use $\bar f$ to replace $f$.

\vskip.5cm
\noindent{\it Proof of 3.2}.
Let $\{p_i'\}$ be $n+1$ points in the standard $n$-sphere $S^n$, and $\{q_i'\}$ the corresponding
anti-podal points satisfying that for $i\neq j$, $d(p_i', p_j')={\frac{\pi}{2}}$. If $d_{GH}(M, S^2)<{\frac{\delta}{3}}$,
from the definition of Gromov-Hausdorff distance, we can then find $p_i, q_i\in M, i=1, 2, \cdots, n+1$, satisfying
\begin{equation*}
|d(p_i, p_j)-{\frac{\pi}{2}}|<\delta,~~ {\text{for}}~ i\neq j,
\end{equation*}
and
\begin{equation*}
|d(p_i, q_i)-\pi|<\delta.
\end{equation*}
Set $g_i(x)=\cos d(p_i, x)$, $i=1, 2, \cdots, n+1$. Then, from Lemma 4, we can find $n+1$ functions $\{f_i, i=1, 2, \cdots, n+1\}$ satisfying
$\Vert f_i\Vert_2\le 1$,
\begin{eqnarray}
\Vert\Delta f_i+nf_i\Vert_2<\psi_1(\delta),
\end{eqnarray}
and
\begin{eqnarray}
\Vert f_i-g_i\Vert_{2,1}<\psi_2(\delta),
\end{eqnarray}
where $\psi_i(\delta)$ satisfy $\lim_{\delta\to 0}\psi_i(\delta)=0$, $i=1, 2$.
By the previous remark, we can assume that $\int_M f_i=0$, $i=1, 2, \cdots, n+1$. If $\{f_i\}$ are linearly
independent, then by applying the minimax principle of eigenvalues of the laplacian
(cf. e.g. \cite{ga-hu-la}, Chapter 4) to the space $\mathcal{H}_0=\{f\in C^\infty(M): \int_Mf=0\}$, the result is obtained.
So, we only need to prove that $\{f_i\}$ are linearly independent as $\delta$ is sufficiently small.

\vskip.5cm
Assume that $\{f_i\}$ are not linearly independent. Then there exist $n+1$ real numbers $a_1, a_2, \cdots, a_{n+1}$
satisfying $\sum_{i=1}^{n+1}a_i^2=1$ and $\sum_{i=1}^{n+1}a_if_i=0$. Set
\begin{eqnarray*}
g(x)=\sum_{i=1}^{n+1}a_ig_i(x)~~~{\text{and}}~~~g'(x')=\sum_{i=1}^{n+1}a_i\cos d(p'_i, x').
\end{eqnarray*}
From (19), we have
\begin{eqnarray}
\Vert g\Vert_{2,1}=\Vert \sum_{i=1}^{n+1}a_ig_i-\sum_{i=1}^{n+1}a_if_i\Vert_{2,1}<\psi_2(\delta).
\end{eqnarray}
On the other hand, by the definition of $g'$, it is an eigenfunction of $S^n$ with eigenvalue being $n$;
so, there exists $q'\in S^n$ such that $g'(x')=\cos d(q', x')$.
It is clear that for $s<{\frac{\pi}{3}}$ and any $q''\in B_s(q')\subset S^n$,
$$
g'(q'')>{\frac 1 2}.
$$
From the choose of $p_i$ and $p_i'$, we also know that if $d_{GH}(x, x')<{\frac{\delta}{3}}$, then
\[
|g(x)-g'(x')|<\delta.
\]
We now choose a $q\in M$ satisfying $d_{GH}(q, q')\le {\frac{\delta}{3}}$,
then for any $\bar q\in B_{\frac s 2}(q)$, we have $g(\bar q)\ge{\frac 1 3}$. But, on the other hand,
the above (20) implies that for sufficiently small $\delta$,
there must exist a $\tilde q\in B_{\frac s 2}(q)$ satisfying $g(\tilde q)<{\frac 1 6}$. This is a contradiction.
The proof is completed.

\vskip.5cm
\noindent

\vskip 1cm
\noindent Yi-Hu Yang: Department of Mathematics, Shanghai Jiao Tong University, Shanghai 200240, China.\\
{\it Email: yangyihu@sjtu.edu.cn}
\\
\\
\noindent Yi Zhang: Department of Mathematics, Tongji University, Shanghai 200092, China.\\
{\it Email: 08zhangyi@tongji.edu.cn}
\end{document}